\documentclass[11pt]{article}
\usepackage{amssymb,amsmath,amsthm,epsfig,times}

\def\C{{\mathbb C}}

\def\N{{\mathbb N}}

\def\R{{\mathbb R}}
\def\*S{{\mathbb S}}

\def\Z{{\mathbb Z}}

\newcommand{\dis}{\displaystyle}

\newcounter{corollaries}

\newtheorem{theorem}{Theorem}
\newtheorem{proposition}{Proposition}[section]
\newtheorem{lemma}[proposition]{Lemma}
\newtheorem{corollary}[corollaries]{Corollary}

\begin{document}

\large

\title{Inversion formulas for complex Radon transform on projective varieties
and boundary value problems for systems of linear PDE.}
\author{Gennadi M. Henkin\footnote{Institut de Mathematiques, Universite Pierre et Marie Curie, 75252, BC247, Paris,
France, and CEMI, Acad. Sc., 117418, Moscow, Russia, henkin@math.jussieu.fr},\hspace{0.05in}
Peter L. Polyakov\footnote{Department of Mathematics, University of Wyoming,
Laramie, WY 82071, USA, polyakov@uwyo.edu}}
\maketitle

\begin{abstract}
Let $G\subset \C P^n$ be a linearly convex compact with smooth boundary, $D={\C}P^n\setminus G$, and let
$D^* \subset \left(\C P^n\right)^*$ be the dual domain. Then for an algebraic,
not necessarily reduced, complete intersection subvariety $V$ of dimension $d$ we construct an explicit inversion
formula for the complex Radon transform $R_V:\ H^{d,d-1}(V\cap D)\to H^{1,0}(D^*)$, and explicit
formulas for solutions of an appropriate boundary value problem for the corresponding system of differential
equations with constant coefficients on $D^*$.
\end{abstract}

\section{Introduction.}\label{Introduction}

Complex Radon-type transforms on complex projective varieties were introduced in different forms and with
different purposes in the works of Fantappie \cite{Fa1}, Martineau \cite{Mar2}, Andreotti, Norguet \cite{AN1,AN2},
Eastwood, Penrose, Wells \cite{Pe, EPW}, Gindikin, Henkin, Polyakov \cite{GH,HP1}, \dots. In a recent paper \cite{HP2}
we have shown that the complex Radon  transform  realizes an isomorphism between the quotient-space
of residual $\bar\partial$-cohomologies $H^{d,d-1}(V\cap D)/H^{d,d-1}(V)$ of algebraic (not necessarily reduced)
$d$-dimensional locally complete intersection $V$ in a linearly concave domain $D$ of ${\C}P^n$ and the space of
holomorphic solutions of the associated homogeneous system of linear
differential equations with constant coefficients in the dual domain $D^*\subset ({\C}P^n)^*$.\\
\indent
In the present paper for an arbitrary algebraic complete intersection $V$ and a smoothly bounded linearly
convex compact $G$ in ${\C}P^n$ we construct an explicit inversion formula for complex Radon transform on $V\cap D$,
where $D={\C}P^n\setminus G$. This inversion formula is based on the explicit formulas for solutions of
appropriate boundary value problems for the associated  with $V$ system of differential equations with constant
coefficients in the dual domain $D^*$. Those formulas are motivated by the  \lq\lq explicit fundamental principle\rq\rq\
of Berndtsson-Passare \cite{BP}.\\
\indent
To formulate the main result of the present paper we introduce the following notations.
Let $(z_0,\ldots,z_n)$ and $(\xi_0,\ldots,\xi_n)$ be the homogeneous
coordinates of points $z\in {\C}P^n$ and $\xi\in ({\C}P^n)^*$. Let
$\langle\xi\cdot z\rangle\stackrel{\rm def}{=}\sum_{k=0}^n\xi_k\cdot z_k$, and let  $\C P^{n-1}_{\xi}$
denote the hyperplane
\begin{equation*}
\C P^{n-1}_{\xi}=\Big\{z\in \C P^n:\ \langle\xi\cdot z\rangle=0\Big\}.
\end{equation*}
Following  \cite{Mar1} and \cite{GH} we call a domain $D\subset \C P^n$ linearly concave,
if there exists a continuous map $D\ni z\to \xi(z)\in \left(\C P^n\right)^*$ satisfying
$$z\in \C P^{n-1}_{\xi(z)}\subset D.$$
A compact $G\subset {\C}P^n$ is called linearly
convex, if the domain $D={\C}P^n\setminus G$ is linearly concave.
The set of hyperplanes, which are contained in the linearly concave domain $D$,
forms the dual domain $D^*\subset ({\C}P^n)^*$. We may assume without loss of generality that the
hyperplane ${\dis \{z\in {\C}P^n:\ z_0=0\} }$ is contained in $D$.\\
\indent
We will denote by $H(D^*,{\cal O}(l))$ the space of holomorphic functions
of homogeneity $l$ on $D^*$. Let $\left\{\tilde P_j\right\}_1^m$ be homogeneous polynomials
of projective coordinates, let $\left\{P_j=\tilde P_j(1,z_1,\dots,z_n)\right\}$ be the corresponding polynomials
of affine coordinates,  and let $V\subset\C P^n$ be the algebraic subvariety
\begin{equation}\label{Variety}
V=\left\{z\in {\C}P^n:\ {\tilde P}_1(z)=\dots ={\tilde P}_m(z)=0\right\}.
\end{equation}
\indent
From \cite{Mar2} we obtain that for $\forall\ g\in H(D^*,{\cal O}(l-1))$ with $l<0$ the solution
$f\in H(D^*,{\cal O}(l))$ of the equation
$$\frac{\partial f}{\partial\xi_0}=g$$
exists and is unique, and therefore the operators
\begin{equation}\label{DOperators}
{\cal D}_j=-\left(\frac{\partial}{\partial\xi_0}\right)^{-1}\frac{\partial}{\partial\xi_j}\
\ \text{for}\ j=1,\dots,n,
\end{equation}
are well defined on the spaces $H(D^*,{\cal O}(l))$ for $l<0$.\\
For a polynomial $R(u)=\sum_{|I|=0}^rR_Iu_1^{i_1}\cdots u_n^{i_n}$ of degree $r$
we denote by $R({\cal D})$ the operator
$$R({\cal D})=\sum_{|I|=0}^rR_I\cdot{\cal D}_1^{i_1}\cdots {\cal D}_n^{i_n}.$$
\indent
We denote by $\{Q^{(k)}\}_{j=1}^m$ the vector-polynomials
$Q^{(k)}(\zeta,z)=\left\{Q^{(k)}_1(\zeta,z),\ldots,Q^{(k)}_n(\zeta,z)\right\}$, such that
$$P_k(\zeta)-P_k(z)=\sum_{j=1}^n(\zeta_j-z_j)Q^{(k)}_j(\zeta,z).$$
\indent
For a linearly convex compact in $\C^n\subset {\C}P^n$
\begin{equation}\label{GCompact}
G=\left\{z\in \C^n:\ \rho(z)\leq 0\right\},
\end{equation}
such that $D={\C}P^n\setminus G$ is a linearly concave domain and $\rho\in C^{\infty}(\C P^n)$, we denote
\begin{align}\label{etaFunctions}
&\eta(\zeta)=\left(\eta_0(\zeta),\eta^{\prime}(\zeta)\right)
=\left(\eta_0(\zeta),\eta_1(\zeta),\ldots,\eta_n(\zeta)\right),\\
&\eta_0(\zeta)=\sum_{j=1}^n\zeta_j\eta_j(\zeta),\
\eta_j(\zeta)=\frac{\partial\rho}{\partial \zeta_j}(\zeta).\nonumber
\end{align}

\begin{theorem}\label{NewTheorem} Let $G$ be a linearly convex compact as in \eqref{GCompact},
$D={\C}P^n\setminus G$, and let $V\subset\C P^n$ be a complete intersection algebraic subvariety as in \eqref{Variety}.\\
\indent
Then any function $g\in {\cal H}({\overline D}^*,{\cal O}(-1))$, satisfying the system of differential equations
\begin{equation}\label{System}
\tilde P_j\left(\frac{\partial}{\partial\xi}\right)g(\xi)=0,\ \text{for}\ j=1,\ldots,m,\ \text{and}\ \xi\in D^*,
\end{equation}
may be represented through its values on the infinitesimal neighborhood of the set
$$\Big\{ \xi\in D^*: \xi=\eta(\zeta)\ \text{for}\ \zeta\in V\cap bG\Big\}$$
by an explicit formula of Cauchy-Fantappie-Leray type:
\begin{multline}\label{gFormula}
g(\xi)=(-1)^{n-m-1}\frac{(n-1)!}{(2\pi i)^n(n-m-1)!}\int_{(\zeta,\mu)\in bG\times\Lambda}
\frac{d\zeta}{\left(\xi_0+\xi^{\prime}\cdot \zeta\right)}\\
\wedge\bar\partial\left(\frac{1}{P_1(\zeta)}\right)
\wedge\cdots\wedge\bar\partial\left(\frac{1}{P_m(\zeta)}\right)
\wedge\omega^{\prime}_0\left(\vartheta(\mu,\zeta,{\cal D})\right)
\Bigg(\frac{\partial^{n-m-1}g}{\partial\eta_0^{n-m-1}}\left(\eta(\zeta)\right)\Bigg),
\end{multline}
where
$$\vartheta(\mu,\zeta,{\cal D})=\sum_{k=1}^m\mu_kQ^{(k)}(\zeta,{\cal D})
+\left(1-\sum_{k=1}^m\mu_k\right)\eta^{\prime}(\zeta)\,$$
$$\omega^{\prime}_0(\vartheta)=\sum_{j=1}^n(-1)^{j-1}\vartheta_j
\wedge_{i\ne j}d\vartheta_i,$$
and the integral in \eqref{gFormula} is understood as
$$\lim_{t\to 0}\int_{T^{\epsilon}_{\left\{{\bf P}\right\}}(t)\times\Lambda}
\frac{d\phi_{\delta}(\zeta)\wedge d\zeta}{\left(\xi_0+\xi^{\prime}\cdot \zeta\right)}
\wedge\frac{\omega^{\prime}_0\left(\vartheta(\mu,\zeta,{\cal D})\right)}
{P_1(\zeta)\ldots P_m(\zeta)}
\Bigg(\frac{\partial^{n-m-1}g}{\partial\eta_0^{n-m-1}}\left(\eta(\zeta)\right)\Bigg)$$
with an arbitrary function $\phi_{\delta}\in {\cal E}_c(\C^n)$ satisfying
\begin{equation}\label{phiFunction}
\phi_{\delta}(\zeta)=\begin{cases}
1\ \mbox{if}\ \rho(\zeta)\leq 0,\vspace{0.1in}\\
0\ \mbox{if}\ \rho(\zeta)> \delta,
\end{cases}
\end{equation}
$$\Lambda=\Big\{\mu\in\R^m_+:\ \sum_{k=1}^m\mu_k\leq 1\Big\},$$
\begin{equation}\label{Ttube}
T^{\epsilon}_{\left\{{\bf P}\right\}}(t)=\Big\{z\in \C^n: |P_k(z)|=\epsilon_k(t)\ \text{for}\ k=1,\dots,m\Big\},
\end{equation}
and $\epsilon(t)=(\epsilon_1(t),\ldots,\epsilon_m(t))$ being an admissible path in the sense of
Coleff-Herrera, i.e. an analytic map $\epsilon:\ [0,1]\to\R_+^m$, satisfying
\begin{equation}\label{admissible}
\lim_{t\to 0}\epsilon_m(t)=0,\ \ \lim_{t\to 0}
\frac{\epsilon_j(t)}{\epsilon_{j+1}^l(t)}=0\ \text{for}\ \forall\ l\in\Z_+.
\end{equation}
\end{theorem}

{\bf Remarks.}
\begin{itemize}
\item
An earlier version of Theorem~\ref{NewTheorem} was proved in \cite{He} for the case 
of the variety $V$ transversally intersecting $bG$, i.e.
$$d\rho\wedge dP_1\wedge\ldots\wedge dP_m\ne 0\ \text{on}\ V\cap bG.$$
\item
Theorem~\ref{NewTheorem} generalizes for the case of general boundary value problems
results of Fantappie \cite{Fa1, Fa2}, Leray \cite{L1, L2}, Rigat \cite{R} on explicit
solutions of the holomorphic Cauchy (or Goursat) problems for systems of linear differential equations
with constant coefficients. Important results on explicit
solutions of nonstandard boundary value problems for two-dimensional
linear integrable PDE were obtained by Fokas \cite{Fo}.
\end{itemize}

A corollary of Theorem~\ref{NewTheorem} presented below is an application of the result of this theorem
to the complex Radon transform. To formulate this corollary we use
definitions from  \cite{HP2}.\\
\indent
A current $f$ in $D$ with support in $V\cap D$ is called a residual current $f\in C^{n-m,n-m-1}(V\cap D)$, if
$$f={\tilde f}\wedge\bar\partial\left(\frac{1}{P_1}\right)\wedge\dots\wedge\bar\partial\left(\frac{1}{P_m}\right),$$
where ${\tilde f}\in {\cal E}^{(n,n-m-1)}$. A residual current $f$ is called $\bar\partial$-closed (denoted
$f\in Z^{n-m,n-m-1}(V\cap D)$ if $\bar\partial{\tilde f}=\sum_{k=1}^mP_k\cdot\Omega_k$ with
$\Omega_k\in {\cal E}^{(n,n-m)}$. We denote by $H^{n-m,n-m-1}(V\cap D)$ the space
$Z^{n-m,n-m-1}(V\cap D)\big/\bar\partial C^{n-m,n-m-2}(V\cap D)$ if $n-m\geq 2$ and
$H^{1,0}(V\cap D)=Z^{1,0}(V\cap D)$ if $n-m=1$.
\begin{corollary}\label{Corollary}
If under the assumptions of Theorem~\ref{NewTheorem}, the coefficient $f_0$ of a closed
holomorphic 1-form $f=\sum_{j=0}^nf_jd\xi_j$ of homogeneity $(-1)$ on $D^*$ satisfies the system
of equations \eqref{System}, then $f$ is the complex Radon transform
$$f(\xi)=R_V[\phi](\xi)=\sum_{j=0}^n\left(\int_{\zeta\in D}
\zeta_j\phi\wedge\bar\partial\frac{1}{\langle\xi \cdot \zeta\rangle}\right)d\xi_j$$
of a residual $\bar\partial$-cohomology class $\phi\in H^{n-m,n-m-1}(V\cap D)$. This cohomology class
corresponds by Serre-Malgrange duality to the functional $\phi^*\in {\cal H}^{\prime}(V\cap G)$,
defined on $\forall\ h\in {\cal H}(V\cap G)$ by the equality
\begin{multline}\label{ResidualFunctional}
\langle \phi^*,h\rangle=(-1)^{m-n-1}\frac{(n-1)!}{(2\pi i)^{n-1}(n-m-1)!}\int_{bG\times\Lambda}h(\zeta)d\zeta
\wedge\bar\partial\left(\frac{1}{P_1(\zeta)}\right)
\wedge\cdots\wedge\bar\partial\left(\frac{1}{P_m(\zeta)}\right)\\
\wedge\omega^{\prime}_0\left(\vartheta(\mu,\zeta,{\cal D})\right)
\Bigg(\frac{\partial^{n-m-1}f_0}{\partial\eta_0^{n-m-1}}\left(\eta(\zeta)\right)\Bigg).
\end{multline}
\end{corollary}

\indent
The proof of Theorem~\ref{NewTheorem} relies on two ingredients: a version of the Martineau type inversion formula
\cite{Mar2,GH} for the Fantappie transform, given here in Proposition~\ref{Martineau}, and an interpolation formula for
holomorphic functions from a complete intersection  subvariety $V\cap G$, not necessarily reduced, to the
linearly convex domain $G\setminus bG$. This interpolation formula \eqref{ExtensionEquality}, proved in
Proposition~\ref{IdealExtension} below, is based on the results of Weil \cite{W}, Leray \cite{L2}, Norguet \cite{N}, and
Coleff-Herrera \cite{CH}.

\section{Cauchy-Leray formula on pseudo-convex complete intersections.}\label{Interpolation}

\indent
In Proposition~\ref{IdealExtension} below we prove a residual interpolation formula in a linearly convex domain,
which can also be considered as the Cauchy-Leray formula for holomorphic functions on complete intersections.
On the one hand the integral term in equality \eqref{ExtensionEquality} of this proposition presents a new interpolation
formula for holomorphic functions in linearly convex domains. On the other hand equality \eqref{ExtensionEquality}
gives a more precise version of the duality theorem of Dickenstein-Sessa and Passare (see \cite{DS,Pa}).

\begin{proposition}\label{IdealExtension} Let $G$ be a linearly convex compact as in \eqref{GCompact},
and let $\{P_k\}_1^m$ be polynomials such that the analytic set
$$V_G=\left\{z\in G:\ P_1(z)=\cdots=P_m(z)=0\right\}$$
is a complete intersection in $G$.
Then for $h\in {\cal H}(G)$ the following formula holds for $z\in G\setminus bG$
\begin{multline}\label{ExtensionEquality}
h(z)=\frac{(n-1)!}{(2\pi i)^n}\Bigg[\lim_{t\to 0}
\int_{T^{\epsilon}_{\left\{{\bf P}\right\}}(t)\times\Lambda}
\frac{h(\zeta)}{\prod_{k=1}^m P_k(\zeta)}d\phi_{\delta}(\zeta)
\wedge\omega^{\prime}_0\Bigg(\sum_{k=1}^m\mu_kQ^{(k)}(\zeta,z)\\
+\left(1-\sum_{k=1}^m\mu_k\right)
\frac{\eta^{\prime}(\zeta)}{\langle\eta^{\prime}(\zeta)\cdot(\zeta-z)\rangle}\Bigg)\wedge d\zeta\Bigg]
+\sum_{k=1}^m h_k(z)\cdot P_k(z),
\end{multline}
where $T^{\epsilon}_{\left\{{\bf P}\right\}}(t)$ is defined in \eqref{Ttube},
$\left\{\epsilon_k(t)\right\}_{k=1}^m$ is an admissible path, function $\phi_{\delta}(\zeta)$ is a function
satisfying \eqref{phiFunction}, $\eta(\zeta)$ is defined in \eqref{etaFunctions}, and $h_k\in H\left(G\right)$.
\end{proposition}
\indent
{\bf Proof.}
We start from the following Weil-Leray-Norguet type integral formula
\begin{multline}\label{StartFormula}
h(z)=\frac{(n-1)!}{(2\pi i)^n}\Bigg[\sum_{0\leq |I|\leq m}\int_{\sigma^{\epsilon}_I(t)\times\Lambda_I}
h(\zeta)\omega^{\prime}_0\Bigg(\sum_{i\in I}\mu_i\frac{Q^{(i)}(\zeta,z)}{P_i(\zeta)-P_i(z)}\\
+\left(1-\sum_{i\in I}\mu_i\right)\frac{\eta^{\prime}(\zeta)}{\langle\eta^{\prime}(\zeta)\cdot(\zeta-z)\rangle}
\Bigg)\wedge\omega(\zeta)\Bigg]
\end{multline}
for a holomorphic function $h$ on the compact
$$U^{\epsilon}(t)=\Big\{z\in \C P^n:\ \rho(z)\leq 0,\
\left\{|P_k(z)|\leq \epsilon_k(t)\right\}_{k=1}^m\Big\},$$
where
$$\sigma_I^{\epsilon}(t)=\Big\{z\in G:\ \rho(z)=0,\
\left\{|P_i(z)|=\epsilon_i(t)\right\}_{i\in I},\ \left\{|P_k(z)|\leq\epsilon_k(t)\right\}_{k\notin I}\Big\},$$
and
$$\Lambda_I=\Big\{\mu\in\R^{|I|}_+:\ \sum_{i\in I}\mu_{i}\leq 1\Big\}.$$
\indent
To transform formula \eqref{StartFormula} into a residue-type formula we assume that function $h$ is defined in
$$G_{\delta}=\left\{z\in \C P^n:\ \rho(z)\leq \delta\right\}$$
for some $\delta>0$, define
$$T^{\epsilon}_I(t)=\Big\{z: 0\leq \rho(z)\leq \delta, \left\{|P_i(z)|=\epsilon_i(t),\right\}_{i\in I},\
\left\{|P_k(z)|\leq\epsilon_k(t),\right\}_{k\notin I}\Big\},$$
and consider the chain
$${\cal C}=\sum_{0\leq |I|\leq m}T^{\epsilon}_I(t)\times\Lambda_I$$
with the boundary
$${\cal B}=\sum_{0\leq |I|\leq m}\left[\sigma^{\epsilon}_I(t)-\sigma^{\epsilon}_I(\delta,t)\right]\times\Lambda_I
+\sum_{0\leq |I|\leq m}T^{\epsilon}_I(t)\times\Gamma_I,$$
where
$$\sigma^{\epsilon}_I(\delta,t)=\Big\{z\in G:\ \rho(z)=\delta,\
\left\{|P_i(z)|=\epsilon_i(t)\right\}_{i\in I},\ \left\{|P_k(z)|\leq\epsilon_k(t)\right\}_{k\notin I}\Big\},$$
and
$$\Gamma_I=\Big\{\mu\in\R^{|I|}_+:\ \sum_{i\in I}\mu_{i}= 1\Big\}.$$
\indent
We consider a function $\phi_{\delta}\in C^{\infty}\left(\C^n\right)$ satisfying \eqref{phiFunction}
and apply the Stokes' formula to the form
\begin{equation*}
h(\zeta)\phi_{\delta}(\zeta)\omega^{\prime}_0\Bigg(\sum_{k=1}^m\mu_k\frac{Q^{(k)}(\zeta,z)}{P_k(\zeta)-P_k(z)}
+\left(1-\sum_{k=1}^m\mu_k\right)
\frac{\eta^{\prime}(\zeta)}{\langle\eta^{\prime}(\zeta)\cdot(\zeta-z)\rangle}
\Bigg)\wedge\omega(\zeta)
\end{equation*}
on the chain ${\cal C}$. Then, using equality $\phi_{\delta}\Big|_{\sigma^{\epsilon}_I(\delta,t)}=0$, we obtain
\begin{multline*}
\sum_{0\leq |I|\leq m}\int_{\sigma^{\epsilon}_I(t)\times\Lambda_I}
h(\zeta)\omega^{\prime}_0\Bigg(\sum_{i\in I}\mu_i\frac{Q^{(i)}(\zeta,z)}{P_i(\zeta)-P_i(z)}
+\left(1-\sum_{i\in I}\mu_i\right)\frac{\eta^{\prime}(\zeta)}{\langle\eta^{\prime}(\zeta)\cdot(\zeta-z)\rangle}
\Bigg)\wedge\omega(\zeta)\\
=-\sum_{0\leq |I|\leq m}\int_{T^{\epsilon}_I(t)\times\Gamma_I}
h(\zeta)\phi_{\delta}(\zeta)\omega^{\prime}_0\Bigg(\sum_{i\in I}\mu_i\frac{Q^{(i)}(\zeta,z)}{P_i(\zeta)-P_i(z)}\\
+\left(1-\sum_{i\in I}\mu_i\right)\frac{\eta^{\prime}(\zeta)}{\langle\eta^{\prime}(\zeta)\cdot(\zeta-z)\rangle}
\Bigg)\wedge\omega(\zeta)\\
+\sum_{0\leq |I|\leq m}\int_{T^{\epsilon}_I(t)\times\Lambda_I}
h(\zeta)d\phi_{\delta}(\zeta)\wedge\omega^{\prime}_0
\Bigg(\sum_{i\in I}\mu_i\frac{Q^{(i)}(\zeta,z)}{P_i(\zeta)-P_i(z)}\\
+\left(1-\sum_{i\in I}\mu_i\right)\frac{\eta^{\prime}(\zeta)}{\langle\eta^{\prime}(\zeta)\cdot(\zeta-z)\rangle}
\Bigg)\wedge\omega(\zeta).
\end{multline*}
\indent
From the dimensional considerations we obtain that 
\begin{multline*}
\int_{T^{\epsilon}_I(t)\times\Gamma_I}
h(\zeta)\phi_{\delta}(\zeta)\omega^{\prime}_0\Bigg(\sum_{i\in I}\mu_i\frac{Q^{(i)}(\zeta,z)}{P_i(\zeta)-P_i(z)}
+\left(1-\sum_{i\in I}\mu_i\right)
\frac{\eta^{\prime}(\zeta)}{\langle\eta^{\prime}(\zeta)\cdot(\zeta-z)\rangle}\Bigg)\wedge\omega(\zeta)\\
=\int_{T^{\epsilon}_I(t)\times\Gamma_I}
h(\zeta)\phi_{\delta}(\zeta)\omega^{\prime}_0
\Bigg(\sum_{i\in I}\mu_i\frac{Q^{(i)}(\zeta,z)}{P_i(\zeta)-P_i(z)}\Bigg)\wedge\omega(\zeta)=0,
\end{multline*}
and therefore the equality above can be rewritten as
\begin{multline}\label{TIntegrals}
\sum_{0\leq |I|\leq m}\int_{\sigma^{\epsilon}_I(t)\times\Lambda_I}
h(\zeta)\omega^{\prime}_0\Bigg(\sum_{i\in I}\mu_i\frac{Q^{(i)}(\zeta,z)}{P_i(\zeta)-P_i(z)}
+\left(1-\sum_{i\in I}\mu_i\right)\frac{\eta^{\prime}(\zeta)}{\langle\eta^{\prime}(\zeta)\cdot(\zeta-z)\rangle}
\Bigg)\wedge\omega(\zeta)\\
=\sum_{0\leq |I|\leq m}\int_{T^{\epsilon}_I(t)\times\Lambda_I}
h(\zeta)d\phi_{\delta}(\zeta)\wedge\omega^{\prime}_0
\Bigg(\sum_{i\in I}\mu_i\frac{Q^{(i)}(\zeta,z)}{P_i(\zeta)-P_i(z)}\\
+\left(1-\sum_{i\in I}\mu_i\right)
\frac{\eta^{\prime}(\zeta)}{\langle\eta^{\prime}(\zeta)\cdot(\zeta-z)\rangle}
\Bigg)\wedge\omega(\zeta).
\end{multline}
\indent
We transform the right-hand side of the last formula for $z\in U^{\epsilon}(t)$ as follows
\begin{multline*}
\sum_{0\leq |I|\leq m}\int_{T^{\epsilon}_I(t)\times\Lambda_I}
h(\zeta)d\phi_{\delta}(\zeta)\wedge\omega^{\prime}_0
\Bigg(\sum_{i\in I}\mu_i\frac{Q^{(i)}(\zeta,z)}{P_i(\zeta)-P_i(z)}\\
+\left(1-\sum_{i\in I}\mu_i\right)
\frac{\eta^{\prime}(\zeta)}{\langle\eta^{\prime}(\zeta)\cdot(\zeta-z)\rangle}
\Bigg)\wedge\omega(\zeta)
\end{multline*}
\begin{multline*}
=\sum_{0\leq |I|\leq m}\int_{T^{\epsilon}_I(t)\times\Lambda_I}
\frac{h(\zeta)}{P_{i_1}(\zeta)}d\phi_{\delta}(\zeta)\wedge\omega^{\prime}_0
\Bigg(\mu_{i_1}Q^{(i_1)}(\zeta,z)
+\sum_{k=2}^{|I|}\mu_{i_k}\frac{Q^{(i_k)}(\zeta,z)}{P_{i_k}(\zeta)-P_{i_k}(z)}\\
+\left(1-\sum_{i\in I}\mu_i\right)
\frac{\eta^{\prime}(\zeta)}{\langle\eta^{\prime}(\zeta)\cdot(\zeta-z)\rangle}\Bigg)\wedge d\zeta\\
+\sum_{0\leq |I|\leq m}\sum_{r=1}^{\infty}\left(P_{i_1}(z)\right)^r
\int_{T^{\epsilon}_I(t)\times\Lambda_I}
\frac{h(\zeta)}{\left(P_{i_1}(\zeta)\right)^{r+1}}d\phi_{\delta}(\zeta)\wedge\omega^{\prime}_0
\Bigg(\mu_{i_1}Q^{(i_1)}(\zeta,z)\\
+\sum_{k=2}^{|I|}\mu_{i_k}\frac{Q^{(i_k)}(\zeta,z)}{P_{i_k}(\zeta)-P_{i_k}(z)}
+\left(1-\sum_{i\in I}\mu_i\right)
\frac{\eta^{\prime}(\zeta)}{\langle\eta^{\prime}(\zeta)\cdot(\zeta-z)\rangle}\Bigg)\wedge d\zeta
\end{multline*}
\begin{multline*}
=\sum_{0\leq |I|\leq m}\int_{T^{\epsilon}_I(t)\times\Lambda_I}
\frac{h(\zeta)}{\prod_{i\in I} P_i(\zeta)}d\phi_{\delta}(\zeta)
\wedge\omega^{\prime}_0\Bigg(\sum_{i\in I}\mu_iQ^{(i)}(\zeta,z)\\
+\left(1-\sum_{i\in I}\mu_i\right)
\frac{\eta^{\prime}(\zeta)}{\langle\eta^{\prime}(\zeta)\cdot(\zeta-z)\rangle}\Bigg)\wedge d\zeta\\
+\sum_{0\leq |I|\leq m}\sum_{r=1}^{\infty}\left(P_{i_1}(z)\right)^r
\int_{T^{\epsilon}_I(t)\times\Lambda_I}
\frac{h(\zeta)}{\left(P_{i_1}(\zeta)\right)^{r+1}}d\phi_{\delta}(\zeta)\wedge\omega^{\prime}_0
\Bigg(\mu_{i_1}Q^{(i_1)}(\zeta,z)\\
+\sum_{k=2}^{|I|}\mu_{i_k}\frac{Q^{(i_k)}(\zeta,z)}{P_{i_k}(\zeta)-P_{i_k}(z)}
+\left(1-\sum_{i\in I}\mu_i\right)
\frac{\eta^{\prime}(\zeta)}{\langle\eta^{\prime}(\zeta)\cdot(\zeta-z)\rangle}\Bigg)\wedge d\zeta
\end{multline*}
$$+\cdots$$
\begin{multline}\label{BigFormula}
+\sum_{0\leq |I|\leq m}\sum_{r=1}^{\infty}\left(P_{i_s}(z)\right)^r
\int_{T^{\epsilon}_I(t)\times\Lambda_I}
\frac{h(\zeta)}{\prod_{k=1}^{s-1}P_{i_k}(\zeta)\left(P_{i_s}(\zeta)\right)^{r+1}}
d\phi_{\delta}(\zeta)\wedge\omega^{\prime}_0\Bigg(\sum_{i\in I}\mu_iQ^{(i)}(\zeta,z)\\
+\left(1-\sum_{i\in I}\mu_i\right)
\frac{\eta^{\prime}(\zeta)}{\langle\eta^{\prime}(\zeta)\cdot(\zeta-z)\rangle}\Bigg)\wedge d\zeta
\end{multline}
where $s=|I|$, and \ \lq\lq\ $\cdots$\ \rq\rq\ stands for the terms of the form
\begin{multline*}
\sum_{0\leq |I|\leq m}\sum_{r=1}^{\infty}\left(P_{i_p}(z)\right)^r
\int_{T^{\epsilon}_I(t)\times\Lambda_I}
\frac{h(\zeta)}{\prod_{k=1}^{p-1}P_{i_k}(\zeta)\left(P_{i_p}(\zeta)\right)^{r+1}}
d\phi_{\delta}(\zeta)\wedge\omega^{\prime}_0\Bigg(\sum_{k=1}^p\mu_{i_k}Q^{(i_k)}(\zeta,z)\\
+\sum_{k=p+1}^{|I|}\mu_{i_k}\frac{Q^{(i_k)}(\zeta,z)}{P_{i_k}(\zeta)-P_{i_k}(z)}
+\left(1-\sum_{i\in I}\mu_i\right)
\frac{\eta^{\prime}(\zeta)}{\langle\eta^{\prime}(\zeta)\cdot(\zeta-z)\rangle}\Bigg)\wedge d\zeta
\end{multline*}
for $1<p<s$.\\
\indent
Denoting then
\begin{multline*}
g_k(z,t)=\frac{(n-1)!}{(2\pi i)^n}\sum_{r=0}^{\infty}\left(P_{k}(z)\right)^r
\sum_{\scriptsize\left\{\begin{array}{ll}
k=i_p\in I\\
1\leq |I|\leq m
\end{array}\right\}}
\int_{T^{\epsilon}_I(t)\times\Lambda_I}
\frac{h(\zeta)}{\prod_{j=1}^{p-1}P_{i_j}(\zeta)\left(P_{k}(\zeta)\right)^{r+1}}d\phi_{\delta}(\zeta)\\
\bigwedge\omega^{\prime}_0\Bigg(\sum_{k=1}^p\mu_{i_k}Q^{(i_k)}(\zeta,z)
+\sum_{k=p+1}^{|I|}\mu_{i_k}\frac{Q^{(i_k)}(\zeta,z)}{P_{i_k}(\zeta)-P_{i_k}(z)}
+\left(1-\sum_{i\in I}\mu_i\right)
\frac{\eta^{\prime}(\zeta)}{\langle\eta^{\prime}(\zeta)\cdot(\zeta-z)\rangle}\Bigg)\wedge d\zeta,
\end{multline*}
and using equalities \eqref{StartFormula}, \eqref{TIntegrals}, and \eqref{BigFormula}
we obtain the following equality for $z\in U^{\epsilon}(t)$
\begin{multline}\label{UFormula}
h(z)=\frac{(n-1)!}{(2\pi i)^n}\Bigg[\sum_{0\leq |I|\leq m}\int_{T^{\epsilon}_I(t)\times\Lambda_I}
\frac{h(\zeta)}{\prod_{i\in I} P_i(\zeta)}d\phi_{\delta}(\zeta)
\wedge\omega^{\prime}_0\Bigg(\sum_{i\in I}\mu_iQ^{(i)}(\zeta,z)\\
+\left(1-\sum_{i\in I}\mu_i\right)
\frac{\eta^{\prime}(\zeta)}{\langle\eta^{\prime}(\zeta)\cdot(\zeta-z)\rangle}\Bigg)\wedge d\zeta\Bigg]
+\sum_{k=1}^mg_k(z,t)\cdot P_k(z).
\end{multline}
\indent
To transform the equality above into equality \eqref{ExtensionEquality} we have to pass to the limit as $t\to 0$
in the right-hand side of \eqref{UFormula}. To prove the existence of limits of the integrals in the right-hand side
of equality above when $t\to 0$ we use the results of Coleff and Herrera. Since all integrals
in \eqref{UFormula} are the integrals of the forms with compact support, those integrals can be reduced to the
integrals of the forms over polydisks. In the proposition below we collect the statements from \cite{CH}, which
are used in the completion of the proof of Proposition~\ref{IdealExtension}.

\begin{proposition}\label{Coleff&Herrera}
Let $D^n=\left\{z\in\C^n:|z_j|<1,\ j=1,\dots,n\right\}$ be
a polydisk in $\C^n$, $\left\{P_k\right\}_{k=1}^m$ - a set of polynomials,
$$V=\left\{z\in D^n: P_1(z)=\cdots=P_m(z)=0\right\}$$
- an algebraic variety of pure dimension $n-m$ such that the restriction to $V$ of the projection
$$\pi:\ D^n\to D^{n-m},$$ defined by the formula $\pi(z_1,\dots,z_n)=(z_{m+1},\dots,z_n)$
is a finite analytic covering, such that the origin is an isolated point in
$\pi^{-1}(0)\cap V$. Let $z^{\prime}=\left(z_1,\dots,z_m\right)$,
and $z^{\prime\prime}=\left(z_{m+1},\dots,z_n\right)$. Then
\begin{enumerate}
\item[(i)]
there exists an analytic function $g$ on $V$ such that for an arbitrary form $\alpha\in {\cal E}_c^{(n,n-m)}\left(D^n\right)$ the following equality holds
\begin{equation}\label{FiberedEquality}
\lim_{t\to 0}\int_{T^{\epsilon}_{\left\{{\bf P}\right\}}(t)}
\frac{\alpha(\zeta)}{\prod_{k=1}^m P_k(\zeta)}
=\lim_{\gamma\to 0}\int_{V\cap\{|g(\zeta)|>\gamma\}}\mbox{res}_{\{{\bf P},\pi\}}[\alpha]\left(\zeta\right),
\end{equation}
where
\begin{equation}\label{TubesandResidues}
\left\{\begin{aligned}
&\mbox{res}_{\{{\bf P},\pi\}}[\alpha]\left(\zeta\right)
=\lim_{t\to 0}\int_{T^{\epsilon}_{\left\{{\bf{\widehat P}}\right\}}(t)}
\frac{{\widehat \alpha}\left(\zeta\right)}{\prod_{k=1}^m {\widehat P}_k(\zeta)},\\
&{\widehat \alpha}\left(\zeta\right)=\alpha\Big|_{\pi^{-1}(\zeta_{m+1},\dots,\zeta_n)},
\hspace{0.05in}{\widehat P}_k=P_k\Big|_{\pi^{-1}(\zeta_{m+1},\dots,\zeta_n)},\hspace{0.05in}
\zeta\in V\cap \pi^{-1}(\zeta_{m+1},\dots,\zeta_n),
\end{aligned}\right.
\end{equation}
and the limit in the left-hand side of \eqref{FiberedEquality} exists,
\item[(ii)]
the limit in the left-hand side of \eqref{FiberedEquality} defines a continuous linear functional on ${\cal E}_c^{(n,n-m)}$,
\item[(iii)]
if $\alpha$ admits a representation $\alpha=f(\zeta)d{\bar \zeta}_{m+1}\wedge d{\bar \zeta}_n\wedge d\zeta$,
then there exist $N\in \N$ and meromorphic functions $\left\{h_I(\zeta)\right\}_{|I|=0}^N$ such that the equality
\begin{equation*}
\mbox{res}_{\{{\bf P},\pi\}}[\alpha]\left(\zeta\right)
=\sum_{|I|=0}^N f_I(\zeta)\cdot h_I(\zeta)
\end{equation*}
holds, where $f_I$ are holomorphic Taylor coefficients of $f$ with respect to $\zeta^{\prime}$.
\item[(iv)]
under conditions of (iii) the following equality holds
\begin{equation}\label{MeromorphicEquality}
\lim_{t\to 0}\int_{T^{\epsilon}_{\left\{{\bf P}\right\}}(t)}
\frac{\alpha(\zeta)}{\prod_{k=1}^m P_k(\zeta)}
=\sum_{|I|=0}^N\lim_{\gamma\to 0}\int_{V\cap\{|g(\zeta)|>\gamma\}} f_I(\zeta)\cdot h_I(\zeta).
\end{equation}
\end{enumerate}
\qed
\end{proposition}

\indent
Using the existence of the limit in  the left-hand side of \eqref{FiberedEquality} we obtain the existence of the limit
$$\lim_{t\to 0}\int_{T^{\epsilon}_I(t)\times\Lambda_I}
\frac{h(\zeta)}{\prod_{i\in I} P_i(\zeta)}d\phi_{\delta}(\zeta)
\wedge\omega^{\prime}_0\Bigg(\sum_{i\in I}\mu_iQ^{(i)}(\zeta,z)
+\left(1-\sum_{i\in I}\mu_i\right)
\frac{\eta^{\prime}(\zeta)}{\langle\eta^{\prime}(\zeta)\cdot(\zeta-z)\rangle}\Bigg)\wedge d\zeta$$
for $I=(1,\dots,m)$.\\
\indent
Also, motivated by equality \eqref{FiberedEquality} we define for $I\subset (1,\dots,m)$ with $|I|=r$ and
$\alpha\in {\cal E}_c^{(n,n-r)}\left(D^n\right)$
\begin{equation}\label{ResidueDefinition}
\lim_{t\to 0}\int_{T^{\epsilon}_I(t)}\frac{\alpha(\zeta)}{\prod_{i\in I} P_i(\zeta)}
=\lim_{\gamma\to 0}\int_{V\cap\{|g(\zeta)|>\gamma\}}\lim_{t\to 0}\int_{{\widehat T}^{\epsilon}_I(\zeta,t)}
\frac{{\widehat \alpha}(w)}{\prod_{i\in I} {\widehat P}_i(w)},
\end{equation}
where we use the notations from \eqref{FiberedEquality}, and additionally
$${\widehat T}^{\epsilon}_I(\zeta,t)=\left\{w\in \pi^{-1}(\zeta_{m+1},\dots,\zeta_n):\
\left\{|{\widehat P}_i(w)|=\epsilon_i(t)\right\}_{i\in I}\right\},\
\left\{|{\widehat P}_k(w)|\leq\epsilon_k(t)\right\}_{k\notin I}.$$
Now, using formula \eqref{ResidueDefinition}, we can pass to the limit as $t\to 0$ in the right-hand side
of \eqref{UFormula} for the integrals from \eqref{UFormula} with $I\neq (1,\dots,m)$. For such integrals
we have the following lemma.

\begin{lemma}\label{ZeroIntegrals} For an arbitrary fixed $z\in G\setminus bG$ the following equality holds
\begin{multline}\label{ZeroEquality}
\lim_{t\to 0}\int_{T^{\epsilon}_I(t)\times\Lambda_I}
\frac{h(\zeta)}{\prod_{i\in I} P_i(\zeta)}d\phi_{\delta}(\zeta)
\wedge\omega^{\prime}_0\Bigg(\sum_{i\in I}\mu_iQ^{(i)}(\zeta,z)\\
+\left(1-\sum_{i\in I}\mu_i\right)
\frac{\eta^{\prime}(\zeta)}{\langle\eta^{\prime}(\zeta)\cdot(\zeta-z)\rangle}\Bigg)\wedge d\zeta=0,
\end{multline}
if  $I\neq (1,\dots,m)$.
\end{lemma}
\indent
{\bf Proof.} To prove equality \eqref{ZeroEquality} we denote
\begin{equation*}
\omega_I(\zeta,z)=\int_{\Lambda_I}\omega^{\prime}_0\Bigg(\sum_{i\in I}\mu_iQ^{(i)}(\zeta,z)
+\left(1-\sum_{i\in I}\mu_i\right)
\frac{\eta^{\prime}(\zeta)}{\langle\eta^{\prime}(\zeta)\cdot(\zeta-z)\rangle}\Bigg)\wedge d\zeta
\end{equation*}
and rewrite the integral in the left-hand side of \eqref{ZeroEquality} as
\begin{equation*}
\lim_{t\to 0}\int_{T^{\epsilon}_I(t)}
\frac{h(\zeta)}{\prod_{i\in I} P_i(\zeta)}d\phi_{\delta}(\zeta)\wedge\omega_I(\zeta,z).
\end{equation*}
\indent
Then using formula \eqref{ResidueDefinition} we rewrite the last limit as
\begin{multline*}
\lim_{t\to 0}\int_{T^{\epsilon}_I(t)}
\frac{h(\zeta)}{\prod_{i\in I} P_i(\zeta)}d\phi_{\delta}(\zeta)\wedge\omega_I(\zeta,z)\\
=\lim_{\gamma\to 0}\int_{V\cap\{|g(\zeta)|>\gamma\}}\lim_{t\to 0}\int_{\scriptsize\left\{\begin{array}{ll}
|{\widehat P}_i(w)|=\epsilon_i(t)\ \text{for}\ i\in I,\\
|{\widehat P}_k(w)|\leq\epsilon_k(t)\ \text{for}\ k\notin I
\end{array}\right\}}
\frac{h(w)}{\prod_{i\in I}{\widehat P}_i(w)}d\phi_{\delta}(w)\wedge\omega_I\left(w,z\right),
\end{multline*}
therefore reducing the proof of the Lemma to the proof of equality
\begin{equation}\label{IsolatedZero}
\lim_{t\to 0}\int_{\scriptsize\left\{\begin{array}{ll}
|{\widehat P}_i(w)|=\epsilon_i(t)\ \text{for}\ i\in I,\\
|{\widehat P}_k(w)|\leq\epsilon_k(t)\ \text{for}\ k\notin I
\end{array}\right\}}
\frac{h(w)}{\prod_{i\in I}{\widehat P}_i(w)}d\phi_{\delta}(w)\wedge\omega_I\left(w,z\right)=0.
\end{equation}
\indent
To prove the last equality we apply the resolution of singularities \cite{Hi} to the isolated point
$$V\cap \pi^{-1}(\zeta^{\prime\prime})=\Bigg\{\zeta^{\prime}\in\C^m:\
{\widehat P}_1(\zeta^{\prime},\zeta^{\prime\prime})=\cdots
={\widehat P}_m(\zeta^{\prime},\zeta^{\prime\prime})=0\Bigg\}$$
in $\C^m(\zeta^{\prime\prime})=\pi^{-1}(\zeta^{\prime\prime})$ for a fixed
$\zeta^{\prime\prime}=(\zeta_{m+1},\dots,\zeta_n)$. Then in a small enough neighborhood
of the origin the lifted variety becomes a normal crossing algebraic variety of the form
$$S=\left\{u\in \C^m:\ u_1^{\alpha_1}=\cdots=u_m^{\alpha_m}=0\right\},$$
and the limit in \eqref{IsolatedZero} becomes
\begin{multline*}
\lim_{t\to 0}\int_{\scriptsize\left\{\begin{array}{ll}
|{\widehat P}_i(w)|=\epsilon_i(t)\ \text{for}\ i\in I,\\
|{\widehat P}_k(w)|\leq\epsilon_k(t)\ \text{for}\ k\notin I
\end{array}\right\}}
\frac{h(w)}{\prod_{i\in I}{\widehat P}_i(w)}d\phi_{\delta}(w)\wedge\omega_I\left(w,z\right)\\
=\lim_{t\to 0}\int_{\scriptsize\left\{\begin{array}{ll}
|u_i^{\alpha_i}|=\epsilon_i(t)\ \text{for}\ i\in I,\\
|u_k^{\alpha_k}|\leq\epsilon_k(t)\ \text{for}\ k\notin I
\end{array}\right\}}
\frac{h^*(u)}{\prod_{i\in I}u_i^{\alpha_i}}d\phi^*_{\delta}(u)\wedge\omega_I^*\left(u,z\right)=0.
\end{multline*}
\qed

\indent
Using Lemma~\ref{ZeroIntegrals} we conclude that in passing to the limit as $t\to 0$ in equality \eqref{UFormula}
the only nonzero may be produced by the integral over $T^{\epsilon}_I(t)\times\Lambda_I$ for $I=(1,\dots,m)$,
i.e. over $T^{\epsilon}_{\left\{{\bf P}\right\}}(t)\times\Lambda$.
The analytic dependence on $z$ of this limit follows from Lemma~\ref{Analyticity} below.

\begin{lemma}\label{Analyticity} Let $D^n$, $V$, $\pi$, and $g$ be the same as in Proposition~\ref{Coleff&Herrera},
and let $T^{\epsilon}_{\left\{{\bf P}\right\}}(t)$ be as in \eqref{Ttube}.\\
\indent
If $F\in {\cal E}^{(n,n-m)}_c\left(D^n\right)$ is a differential form with respect to variables $\zeta$, with coefficients
infinitely differentiable with respect to both variables $\zeta$ and $z$, and holomorphic with respect to variables $z$,
then
\begin{equation*}
R(z)=\lim_{t\to 0}\int_{T^{\epsilon}_{\left\{{\bf P}\right\}}(t)}\frac{F(\zeta,z)}{\prod_{k=1}^m P_k(\zeta)},
\end{equation*}
is a holomorphic function.
\end{lemma}
\indent
{\bf Proof.} Without loss of generality we may assume that
$F(\zeta,z)=f(\zeta,z)d{\bar \zeta}_{m+1}\wedge d{\bar \zeta}_n\wedge d\zeta$.
Then, following \cite{CH}, we consider the Taylor series of $f$ at $\zeta\in V$ with respect to $\zeta^{\prime}$
\begin{equation*}
f(w,z)\Big|_{\pi^{-1}(\zeta^{\prime\prime})}=\sum_{|I|+|J|=0}^{\infty}f_{I,J}(\zeta,z)\cdot
\left(w^{\prime}-\zeta^{\prime}\right)^I\cdot \left(\overline{w}^{\prime}-\overline{\zeta}^{\prime}\right)^J,
\end{equation*}
and using equality \eqref{MeromorphicEquality} obtain the existence of $N\in \N$ and of meromorphic functions
$\left\{h_I(\zeta)\right\}_{|I|=0}^N$ such that the equality
\begin{equation}\label{MeromorphicSum}
\lim_{t\to 0}\int_{T^{\epsilon}_{\left\{{\bf P}\right\}}(t)}\frac{F(\zeta,z)}{\prod_{k=1}^m P_k(\zeta)}
=\sum_{|I|=0}^N\lim_{\gamma\to 0}\int_{V\cap\{|g(\zeta)|>\gamma\}} f_I(\zeta,z)\cdot h_I(\zeta)
\end{equation}
holds.\\
\indent
If $f(\zeta,z)$ is a polynomial with respect to $z$, then the left-hand side of \eqref{MeromorphicSum} is a polynomial
as well. For an arbitrary $f(\zeta,z)\in {\cal E}_c\left(D^n\right)$ analytically depending on $z$ we approximate
it by polynomials, and then use the continuity of a residual current as a functional on ${\cal E}_c^{(n,n-m)}$, which
follows from (ii) in Proposition~\ref{Coleff&Herrera}.
\qed\\

\indent
Continuing with the proof of Proposition~\ref{IdealExtension} we obtain from Lemmas~\ref{ZeroIntegrals} and
\ref{Analyticity} the following equality
\begin{multline*}
\sum_{0\leq |I|\leq m}\lim_{t\to 0}\int_{T^{\epsilon}_I(t)\times\Lambda_I}
\frac{h(\zeta)}{\prod_{i\in I} P_i(\zeta)}d\phi_{\delta}(\zeta)
\wedge\omega^{\prime}_0\Bigg(\sum_{i\in I}\mu_iQ^{(i)}(\zeta,z)\\
+\left(1-\sum_{i\in I}\mu_i\right)
\frac{\eta^{\prime}(\zeta)}{\langle\eta^{\prime}(\zeta)\cdot(\zeta-z)\rangle}\Bigg)\wedge d\zeta\\
=\lim_{t\to 0}\int_{T^{\epsilon}_{\left\{{\bf P}\right\}}(t)\times\Lambda}
\frac{h(\zeta)}{\prod_{i=1}^m P_i(\zeta)}d\phi_{\delta}(\zeta)
\wedge\omega^{\prime}_0\Bigg(\sum_{i=1}^m\mu_iQ^{(i)}(\zeta,z)\\
+\left(1-\sum_{i=1}^m\mu_i\right)
\frac{\eta^{\prime}(\zeta)}{\langle\eta^{\prime}(\zeta)\cdot(\zeta-z)\rangle}\Bigg)\wedge d\zeta,
\end{multline*}
with the right-hand side being holomorphic with respect to $z$.\\
\indent
To prove the existence of coefficients $h_k\in H\left(G\right)$ in \eqref{ExtensionEquality},
and therefore to complete the proof of Proposition~\ref{IdealExtension} we notice that the functions
\begin{multline*}
u_t(z)=h(z)-\frac{(n-1)!}{(2\pi i)^n}\Bigg[\sum_{0\leq |I|\leq m}\int_{T^{\epsilon}_I(t)\times\Lambda_I}
\frac{h(\zeta)}{\prod_{i\in I} P_i(\zeta)}d\phi_{\delta}(\zeta)\\
\wedge\omega^{\prime}_0\Bigg(\sum_{i\in I}\mu_iQ^{(i)}(\zeta,z)
+\left(1-\sum_{i\in I}\mu_i\right)
\frac{\eta^{\prime}(\zeta)}{\langle\eta^{\prime}(\zeta)\cdot(\zeta-z)\rangle}\Bigg)\wedge d\zeta\Bigg]
\end{multline*}
form a family of holomorphic functions on the interior of $G$ depending on $t$ and converging to the
holomorphic function
\begin{multline*}
u(z)=h(z)-\frac{(n-1)!}{(2\pi i)^n}\Bigg(\lim_{t\to 0}\int_{T^{\epsilon}_J(t)\times\Lambda_J}
\frac{h(\zeta)}{\prod_{i=1}^m P_i(\zeta)}d\phi_{\delta}(\zeta)\\
\wedge\omega^{\prime}_0\left(\sum_{i=1}^m\mu_iQ^{(i)}(\zeta,z)
+\left(1-\sum_{i=1}^m\mu_i\right)
\frac{\eta^{\prime}(\zeta)}{\langle\eta^{\prime}(\zeta)\cdot(\zeta-z)\rangle}\right)\wedge d\zeta\Bigg)
\end{multline*}
on the interior of $G$. Since for each $t$ the function $u_t$ defines a section
of the sheaf of ideals, defined by the functions $P_1,\dots,P_m$ on $G$
from H. Cartan's Theorems (A) and (B) in \cite{Ca} we obtain that the limit function
$u=\lim_{t\to 0}u_t$ admits a representation on the interior of $G$
$$u(z)=\sum_{k=1}^m h_k(z)\cdot P_k(z)$$
with $h_k\in H(G)$.\qed\\
\section{Proof of Theorem~\ref{NewTheorem}.}\label{Proof}

\indent
Before proving Theorem~\ref{NewTheorem} we present a version of the Martineau's (see \cite{Mar2})
inversion formula for the Fantappi\'e transform from \cite{GH}, which is used in the proof.\\
\indent
For $f\in {\cal H}\left(D^*\right)$, following \cite{Mar2} and \cite{GH}, we consider the analytic functional
$\mu^f$ on the space $H(G)$ defined by the formula
\begin{equation}\label{mu_f}
\mu^f(h)=\int_{bG_{-\nu}}h\cdot \Omega_f,
\end{equation}
where
$$\Omega_f(z)=\frac{(-1)}{(2\pi i)^n}\ \frac{\partial^n f}{\partial\eta_0^n}(\eta(z))
\omega^{\prime}\left(\eta(z)\right)\bigwedge_{j=1}^nd\left(\frac{z_j}{z_0}\right),$$
$G_{-\nu}=\left\{z\in\C P^n: \rho(z)\leq -\nu\right\}$,
and a map $\eta:bG_{-\nu}\to \left(\C P^n\right)^*$ satisfies
$\langle\eta(z)\cdot z\rangle=0$ for $z\in bG_{-\nu}$.\\
\indent
The indicatrice of Fantappi\'e of the functional $\mu^f$ is a holomorphic $1$-form on $D^*$ defined by the formula
\begin{multline*}
{\cal F}\mu^f=\sum_{k=0}^n\mu^f\left(\frac{z_k}{\langle\xi\cdot z\rangle}\right)d\xi_k\\
=\frac{(-1)}{(2\pi i)^n}\sum_{k=0}^n\left(\int_{bG_{-\nu}}
\left(\frac{z_k}{\langle\xi\cdot z\rangle}\right)
\frac{\partial^n f}{\partial\eta_0^n}(\eta(z))\omega^{\prime}\left(\eta(z)\right)
\bigwedge_{j=1}^nd\left(\frac{z_j}{z_0}\right)\right)d\xi_k.
\end{multline*}
\indent
The most important application of the indicatrice of Fantappi\'e of $\mu^f$ is the inversion formula
described in the proposition below.

\begin{proposition}\label{Martineau}(Martineau type inversion formula.
\cite{Mar2}, \cite{GH}.)\ Let $D\subset \C P^n$ be a linearly concave domain such that
$D^*\subset \left\{\xi_0\neq 0\right\}$, and let $f\in {\cal H}(D^*)$. Then the following equality holds:
\begin{equation}\label{MartineauFormula}
{\cal F}\mu^f(\xi)=df(\xi),
\end{equation}
or
$$\frac{(-1)}{(2\pi i)^n}\int_{z\in bG_{-\nu}}
\frac{z_k}{\langle\xi\cdot z\rangle}
\frac{\partial^{n}f}{\partial\eta_0^{n}}(\eta(z))
\omega^{\prime}\left(\eta(z)\right)\bigwedge_{j=1}^nd\left(\frac{z_j}{z_0}\right)
=\frac{\partial f}{\partial\xi_k}(\xi)$$
for $k=0,\dots,n$, and $\xi\in D^*$.\\
\indent
Moreover, for $g\in {\cal H}(D^*,{\cal O}(-1))$ we have the following equality
\begin{equation}\label{Inverse}
g(\xi)=\frac{(-1)}{(2\pi i)^n}\int_{bG_{-\nu}}
\frac{\partial^{n-1}g}{\partial\eta_0^{n-1}}(\eta(u))
\frac{\omega^{\prime}\left(\eta(u)\right)\wedge du}{\left(\xi_0+\xi^{\prime}
\cdot u\right)},
\end{equation}
where $\xi\in D^*$ and
$$u_j=\frac{z_j}{z_0}\ \mbox{for}\ j=1,\dots,n.$$
\end{proposition}
\qed

\indent
To prove Theorem~\ref{NewTheorem} we consider $g\in {\cal H}(D^*,{\cal O}(-1))$ satisfying
the system of equations \eqref{System} and using equality \eqref{Inverse} obtain the equality
\begin{multline*}
\frac{(-1)^{1+\deg{P}_k}}{(2\pi i)^n}\left(\deg{P}_k\right)!
\int_{bG_{-\nu}}\frac{P_k(u)}{\left(\xi_0+\xi^{\prime}\cdot u\right)^{1+\deg{P}_k}}
\cdot\frac{\partial^{n-1}g}{\partial\eta_0^{n-1}}(\eta(u))
\omega^{\prime}\left(\eta(u)\right)\wedge du\\
=P_k\left(\frac{\partial}{\partial\xi}\right)\left[g\right](\xi)=0.
\end{multline*}
Then, using the Cauchy-Leray formula \cite{L2} we obtain the density of the set of functions
$$\left\{\frac{1}{\left(\xi_0+\xi^{\prime}\cdot u\right)^{1+\deg{P}_k}}\right\}_{\xi\in D^*}$$
in the space $H(G)$, and therefore the equality
\begin{equation}\label{ZeroOnIdeal}
\int_{bG_{-\nu}}f(u)\cdot P_k(u)
\cdot\frac{\partial^{n-1}g}{\partial\eta_0^{n-1}}(\eta(u))
\omega^{\prime}\left(\eta(u)\right)\wedge du=0
\end{equation}
for an arbitrary $f\in H(G)$.\\
\indent
Using notation \eqref{etaFunctions} for $\eta_0(w)=\left\langle\eta^{\prime}(w)\cdot w\right\rangle$ and applying
Proposition~\ref{IdealExtension} we consider the function
\begin{multline}\label{HFunction}
H_V(\xi,u)=\frac{(n-1)!}{(2\pi i)^n}\lim_{t\to 0}
\int_{T^{\epsilon}_{\left\{{\bf P}\right\}}(t)\times\Lambda}
\frac{d\phi_{\delta}(w)}{\prod_{k=1}^m P_k(w)\cdot\left(\xi_0+\xi^{\prime}\cdot w\right)}\\
\wedge\omega^{\prime}_0\left(\sum_{k=1}^m\mu_k Q^{(k)}(w,u)
+\left(1-\sum_{k=1}^m\mu_k\right)
\frac{\eta^{\prime}(w)}{\langle\eta^{\prime}(w)\cdot (w-u)\rangle}\right)\wedge dw\\
=\frac{(n-1)!}{(2\pi i)^n}\lim_{t\to 0}\int_{T^{\epsilon}_{\left\{{\bf P}\right\}}(t)\times\Lambda}
\frac{d\phi_{\delta}(w)}{\prod_{k=1}^m P_k(w)\cdot\left(\xi_0+\xi^{\prime}\cdot w\right)}\\
\wedge\omega^{\prime}_0\left(\sum_{k=1}^m\mu_k Q^{(k)}(w,u)
+\left(1-\sum_{k=1}^m\mu_k\right)
\frac{\eta^{\prime}(w)}{\left(\eta_0(w)-\langle\eta^{\prime}(w)\cdot u\rangle\right)}\right)\wedge dw,
\end{multline}
satisfying the equality
$$H_V(\xi,u)=\frac{1}{\left(\xi_0+\xi^{\prime}\cdot u\right)}
+\sum_{k=1}^mh_k(\xi,u)\cdot P_k(u)$$
for $u\in G$.\\
\indent
Using the equality above and equality (\ref{ZeroOnIdeal}) in equality (\ref{Inverse})
we obtain the following equality
\begin{multline*}
g(\xi)=\frac{(-1)}{(2\pi i)^n}\int_{bG_{-\nu}}
\frac{\partial^{n-1}g}{\partial\eta_0^{n-1}}(\eta(u))
\frac{\omega^{\prime}\left(\eta(u)\right)\wedge du}{\left(\xi_0+\xi^{\prime}\cdot u\right)}\\
=\frac{(-1)}{(2\pi i)^n}\int_{bG_{-\nu}}
\frac{\partial^{n-1}g}{\partial\eta_0^{n-1}}(\eta(u))
H_V(\xi,u)\omega^{\prime}\left(\eta(u)\right)\wedge du,
\end{multline*}
which after the substitution of expression (\ref{HFunction}) and the change
of the order of integration becomes
\begin{multline}\label{gFormula1}
g(\xi)=\frac{(n-1)!}{(2\pi i)^n}\lim_{t\to 0}
\int_{T^{\epsilon}_{\left\{{\bf P}\right\}}(t)\times\Lambda}
\frac{d\phi_{\delta}(w)\wedge dw}{\prod_{k=1}^m P_k(w)\cdot\left(\xi_0+\xi^{\prime}\cdot w\right)}\\
\times\frac{(-1)}{(2\pi i)^n}\int_{bG_{-\nu}}
\frac{\partial^{n-1}g}{\partial\eta_0^{n-1}}(\eta(u))
\wedge\omega^{\prime}_0\left(\sum_{k=1}^m\mu_k Q^{(k)}(w,u)\right.\\
\left.+\left(1-\sum_{k=1}^m\mu_k\right)
\frac{\eta^{\prime}(w)}{\left(\eta_0(w)-\langle\eta^{\prime}(w)\cdot u\rangle\right)}\right)
\wedge\omega^{\prime}\left(\eta(u)\right)\wedge du.
\end{multline}
\indent
To transform formula \eqref{gFormula1} we notice that the operators ${\cal D}_j$ defined in \eqref{DOperators},
satisfy the following condition
\begin{multline}\label{DCondition}
{\cal D}_i\left(\frac{1}{\eta_0(w)-\sum_{j=1}^n\eta_j(w)u_j}\right)
=-\left(\frac{\partial}{\partial\eta_0}\right)^{-1}
\left(\frac{u_i}{\Big(\eta_0(w)-\sum_{j=1}^n\eta_j(w)u_j\Big)^2}\right)\\
=\frac{u_i}{\eta_0(w)-\sum_{j=1}^n\eta_j(w)u_j}.
\end{multline}
Then, using equality \eqref{DCondition} we obtain that for a polynomial $Q^{(k)}_j(w,u)$
the differential operator $Q^{(k)}_j(w,{\cal D})$ satisfies the following property
\begin{equation*}
Q^{(k)}_j(w,{\cal D})\left(\frac{1}{\eta_0(w)-\langle\eta^{\prime}(w)\cdot u\rangle}\right)
=\frac{Q^{(k)}_j(w,u)}{\eta_0(w)-\langle\eta^{\prime}(w)\cdot u\rangle}.
\end{equation*}
\indent
Using the equality above we rewrite equality \eqref{gFormula1} as
\begin{multline*}
g(\xi)=\frac{(n-1)!}{(2\pi i)^n}\lim_{t\to 0}
\int_{T^{\epsilon}_{\left\{{\bf P}\right\}}(t)\times\Lambda}
\frac{d\phi_{\delta}(w)\wedge dw}{\prod_{k=1}^m P_k(w)\cdot\left(\xi_0+\xi^{\prime}\cdot w\right)}\\
\wedge\omega^{\prime}_0\Bigg(\sum_{k=1}^m\mu_k Q^{(k)}(w,{\cal D})
+\left(1-\sum_{k=1}^m\mu_k\right)\eta^{\prime}(w)\Bigg)\\
\frac{(-1)}{(2\pi i)^n}\int_{bG_{-\nu}}
\frac{\partial^{n-1}g}{\partial\eta_0^{n-1}}(\eta(u))
\frac{\omega^{\prime}\left(\eta(u)\right)\wedge du}
{\left(\eta_0(w)-\langle\eta^{\prime}(w)\cdot u\rangle\right)^{n-m}}
\end{multline*}
\begin{multline*}
=\frac{(n-1)!}{(2\pi i)^n}\lim_{t\to 0}
\int_{T^{\epsilon}_{\left\{{\bf P}\right\}}(t)\times\Lambda}
\frac{d\phi_{\delta}(w)\wedge dw}{\prod_{k=1}^m P_k(w)\cdot\left(\xi_0+\xi^{\prime}\cdot w\right)}\\
\wedge\omega^{\prime}_0\left(\sum_{k=1}^m\mu_k Q^{(k)}(w,{\cal D})
+\left(1-\sum_{k=1}^m\mu_k\right)\eta^{\prime}(w)\right)\\
\frac{\partial^{n-m-1}}{\partial\eta_0^{n-m-1}}
\frac{(-1)^{n-m-1}}{(m-n-1)!}\left[\frac{(-1)}{(2\pi i)^n}\int_{bG_{-\nu}}
\frac{\partial^{n-1}g}{\partial\eta_0^{n-1}}(\eta(u))
\frac{\omega^{\prime}\left(\eta(u)\right)\wedge du}
{\left(\eta_0(w)-\langle\eta^{\prime}(w)\cdot u\rangle\right)}\right]
\end{multline*}
\begin{multline*}
=(-1)^{m-n-1}\frac{(n-1)!}{(2\pi i)^n(n-m-1)!}\lim_{t\to 0}
\int_{T^{\epsilon}_{\left\{{\bf P}\right\}}(t)\times\Lambda}
\frac{d\phi_{\delta}(w)\wedge dw}{\prod_{k=1}^m P_k(w)\cdot\left(\xi_0+\xi^{\prime}\cdot w\right)}\\
\wedge\omega^{\prime}_0\left(\sum_{k=1}^m\mu_k Q^{(k)}(w,{\cal D})
+\left(1-\sum_{k=1}^m\mu_k\right)\eta^{\prime}(w)\right)
\left(\frac{\partial^{n-m-1}g}{\partial\eta_0^{n-m-1}}\left(\eta(w)\right)\right),
\end{multline*}
where in the last equality we have used equality \eqref{Inverse}.\qed\\

{\bf Proof of Corollary~\ref{Corollary}.}

\indent
Let $f_0$ be the coefficient of a closed 1-form $f=\sum_{j=0}^nf_jd\xi_j$ on $D^*$ satisfying
the system of equations \eqref{System}. Then from equality \eqref{gFormula} in Theorem~\ref{NewTheorem} 
for the functional $\phi^*$ defined in \eqref{ResidualFunctional} we obtain the equality
$$f_0=\frac{1}{2\pi i}\left\langle \phi^*,\frac{1}{\xi_0+\xi^{\prime}w}\right\rangle.$$
\indent
On the other hand, using the linear convexity of $D^*$, we can find $g\in H(D^*)$ such that
$f=dg$, and, in particular, $f_0={\partial g}/{\partial\xi_0}$. Since the function
$g$ has homogeneity $0$, the following equality holds
$$\xi_0\frac{\partial g}{\partial\xi_0}=-\sum_{j=1}^n\xi_j\frac{\partial g}{\partial\xi_j},$$
which leads to equality $\langle \phi^*,1\rangle=f_0(1,0,\ldots,0)=0$. From the closedness of the form $f$
we obtain the following equality
$$\frac{1}{2\pi i}\sum_{j=0}^n\left\langle \phi^*,\frac{w_j}{{\xi_0+\xi^{\prime}w}}\right\rangle d\xi_j
=\sum_{j=0}^nf_jd\xi_j=f,\ w\in G,\ w_0=1,\ \xi\in D^*.$$
\indent
Since a complete intersection $V$ in ${\C}P^n$ is connected (see \cite{Ha}, $\S$ III.5), Theorem 2 of \cite{HP2}
implies the existence of a residual cohomology class $\phi\in H^{n-m,n-m-1}(V\cap D)$ such that
$$R_V[\phi](\xi)=f(\xi)=\frac{1}{2\pi i}
\sum_{j=0}^n\left\langle \phi^*,\frac{w_j}{{\xi_0+\xi^{\prime}w}}\right\rangle d\xi_j.$$
\indent
A representative of this cohomology class can be found explicitly. To find such a representative
we consider $\phi^*$ as the $(n-m,n-m)$-current with support in $V\cap G$. Condition $\langle \phi^*,1\rangle=0$
and the connectedness of $V$ imply by Serre-Malgrange duality (see \cite{S}, \cite{Mal}) that
there exists a current ${\hat\phi}$ of bidegree  $(n-m,n-m-1)$ on $V$, such that $\bar\partial{\hat\phi}=\phi^*$ on $V$.
So the sought cohomology class on $V\cap D$ can be defined as $\phi={\hat\phi}\big|_{V\cap D}$.
\qed


\begin{thebibliography}{WIDT}
\small
\bibitem[AN1]{AN1} A. Andreotti, F. Norguet, La convexit\'e holomorphe dans
l'espace analytique des cycles d'une vari\'et\'e alg\'ebrique, Ann. Scuola
Norm. Sup. Pisa 21 (1967), 31-82.
\bibitem[AN2]{AN2} A. Andreotti, F. Norguet, Cycles of algebraic manifolds
and $\partial\bar\partial$-cohomology, Ann. Scuola Norm. Sup. Pisa 25 (1971), 59-114.
\bibitem[BP]{BP} B. Berndtsson, M. Passare, Integral formulas and an explicit
version of the fundamental principle, J. of Functional Analysis, 84 (1989), 358-372.
\bibitem[Ca]{Ca} H. Cartan, S\'eminaire E.N.S., 1951-1952, \'Ecole Normale
Sup\'erieure, Paris.
\bibitem[CH]{CH} N.R. Coleff, M.E. Herrera, Les Courants R\'esiduels Associ\'es \`a une
FormeM\'eromorphe, Lecture Notes in Mathematics, 633, Springer Verlag, New York, 1978.
\bibitem[DS]{DS} A. Dickenstein, C. Sessa, Canonical representatives in moderate cohomology,
Invent. Math. 80 (1985), 417-434.
\bibitem[EPW]{EPW} M. Eastwood, R. Penrose, R.O. Wells, Cohomology and
massless fields, Comm. Math. Phys. 78 (1980/1981), 305-351.
\bibitem[Fa1]{Fa1} L. Fantappie, L'indicatrice proiettiva dei funzionali
ei prodotti funzionali proiettivi, Annali di Mat., 13 (1943), 1-100.
\bibitem[Fa2]{Fa2} L. Fantappie, Sur les m\'ethodes nouvelles d'int\'egration
 des \'equations aux d\'eriv\'ees partielles au moyen des fonctionnelles
analytiques, La th\'eorie des \'equations aux d\'eriv\'ees partielles 81, Coll.Int. CNRS, Nancy, 1956.
\bibitem[Fo]{Fo} A.S. Fokas, On the integrability of  linear and
nonlinear partial differential equations, J. of Math. Physics 41(6) (2000),
4188-4237.
\bibitem[GH]{GH} S. Gindikin, G. Henkin, Integral geometry for $\bar\partial$-cohomology
in $q$-linearly concave domains in $\C P^n$, Funct. Anal. Appl. 12 (1979), 247-261.
\bibitem[Ha]{Ha} R. Hartshorne, Algebraic Geometry, Springer-Verlag, New York, 1977.
\bibitem[He]{He} G. M. Henkin, The Abel-Radon transform and several complex variables,
Ann. of Math. Studies, 137 (1995), 223-275, Preprint University Paris VI, 1993.
\bibitem[HP1]{HP1} G.M. Khenkin=Henkin, P.L. Polyakov, Homotopy formulas for the
$\bar\partial$ - operator on $\C P^n$ and the Radon-Penrose transform, Math. USSR Izvestiya,
28:3 (1987), 555-587.
\bibitem[HP2]{HP2} G.M. Henkin, P.L. Polyakov, Residual $\bar\partial$-cohomology and the complex
Radon transform on subvarieties of $\C P^n$, arXiv:1012.4438v1 (2010).
\bibitem[Hi]{Hi} H. Hironaka, The resolution of singularities of an algebraic variety over a field
of characteristic zero, Ann. Math. 19 (1964), 109 - 326.
\bibitem[L1]{L1} J. Leray, Hyperbolic differential equations, The Inst. for Advanced Study, Princeton, 1953.
\bibitem[L2]{L2} J. Leray, Le calcul differ\'entiel et int\'egral sur une vari\'et\'e
analytique complexe, Bull. Soc. Math. France 87 (1959), 81-180.
\bibitem[Mar1]{Mar1} A. Martineau, Indicatrices des fonctionelles analytiques et inversion
de la transformation de Fourier-Borel par la transformation de Laplace, C.R. Acad. Sci. Paris 255
(1962), 1845-1847, 2888-2890.
\bibitem[Mar2]{Mar2} A. Martineau, Equations diff\'erentiels d'ordre infini, Bull. Soc. Math.
France, 95 (1967), 109-154.
\bibitem[Mal]{Mal} B. Malgrange, Syst\`emes diff\'erentiels \`a coefficients
constants. S\'eminaire Bourbaki 15(246) (1962-1963).
\bibitem[N]{N} F. Norguet, Probl\`ems sur les formes diff\'erentielles et des
courants, Ann. Inst. Fourier 11 (1961).
\bibitem[Pa]{Pa} M. Passare, Residues, currents, and their relation to ideals of holomorphic functions,
Math. Stand. 62 (1988), 75-152.
\bibitem[Pe]{Pe} R. Penrose, Massless fields and sheaf cohomology, Twistor
Newsletter 5, July 1977, Oxford.
\bibitem[R]{R} S. Rigat, Application of the fundamental principle to
complex Cauchy problem, Ark. Mat. 38(2) (2000), 355-380.
\bibitem[S]{S} J.P. Serre, Un th\'eor\`eme de dualit\'e, Comm. Math. Helv. 29 (1955), 9-26.
\bibitem[W]{W} A. Weil, L'int\'egral de Cauchy et les fonctions de plusieurs
variables, Math. Annalen 111 (1935), 178-182.
\end{thebibliography}
\end{document}